\documentclass{article}

\usepackage{amssymb}

\usepackage{amsthm,amsmath}

\newtheorem{theorem}{Theorem}
\newtheorem{lemma}[theorem]{Lemma}

\newtheorem{conjecture}[theorem]{Conjecture}

\newtheorem{proposition}[theorem]{Proposition}

\newtheorem{question}[theorem]{Question}

\newcommand{\Z}{\mathbf{Z}}
\newcommand{\Q}{\mathbf{Q}}
\newcommand{\F}{\mathbf{F}}
\DeclareMathOperator{\SL}{SL}
\newcommand{\sltz}{\SL_2(\Z)}

\begin{document}

\title{On mod~$p$ modular representations which are defined over~$\F_p$}
\author{L. J. P. Kilford\footnote{Keywords: Modular forms, Hecke polynomials} \footnote{Mathematics Subject Classification: 11F11, 11F33.}}

\maketitle

\begin{abstract}
In this paper, we will use techniques of Conrey, Farmer and Wallace to find spaces of modular forms~$S_k(\Gamma_0(N))$ where all of the eigenspaces have Hecke eigenvalues defined over~$\F_p$, and give a heuristic that these are all such spaces.
\end{abstract}

\section{Introduction}
Let~$N$ and~$k$ be positive integers and let~$l$ and~$p$ be primes. We will let~$S_k(\Gamma_0(N),\chi)$ be the space of holomorphic cusp forms of integer weight~$k$ for the congruence subgroup~$\Gamma_0(N)$, and we will define~$T_{l,k}^{N,\chi}$ to be the characteristic polynomial of the Hecke operator~$T_l$ acting on~$S_k(\Gamma_0(N),\chi)$. We will call this polynomial the \emph{Hecke polynomial}.

For modular forms in characteristic~0, there is a well-known conjecture (Maeda's conjecture) that says that the characteristic polynomials of the Hecke operators acting on modular forms for~$\sltz$ are irreducible:
\begin{conjecture}[Maeda's Conjecture]
Let~$k$ be a positive integer, and let~$l$ be a prime number. The Hecke polynomial of~$T^{1,1}_{l,k}$ is irreducible, with Galois group~$S_n$, where~$n$ is the dimension of~$S_k(\sltz)$ as a complex vector space.
\end{conjecture}
Some progress has been made towards this conjecture; methods introduced in~\cite{buzzard-t2} prove that certain Hecke polynomials are irreducible and have full Galois group, and results such as those in~\cite{conrey-farmer-wallace}, \cite{baba-murty} and~\cite{ahlgren-new} show that if a certain~$T_l$ is irreducible then other~$T_r$ must be irreducible also.

In the paper~\cite{conrey-farmer-wallace}, a method for showing that certain of the~$T_{l,k}^{1,1}$ are irreducible over~$\Q$ is given; this involves a careful study of how level~1 Hecke polynomials factor modulo certain primes which uses the fact that the factorizations of the Hecke polynomials are periodic. In~\cite{james-and-ono} and~\cite{baba-murty}, the authors first show that certain of the~$T^{1,1}_{2,k}$ are irreducible and then using effective versions of the Chebotarev density theorem to show that a positive proportion of the~$T^{1,1}_{l,k}$ are irreducible, and in~\cite{ahlgren-new}, plausible hypotheses on the non-vanishing of certain coefficients of modular forms are used to prove irreducibility of the Hecke polynomials (when the weight is~12, for instance, these hypotheses are exactly Lehmer's conjecture on the Ramanujan $\tau$-function).%the structure of the space of modular forms of level~1 to show that  of James and Ono~\cite{james-and-ono} to show that the other~$T^{1,1}_{l,k}$ are irreducible. %this involves first showing that such that~$T^{1,1}_{2,k}$ is irreducible and then using the structure of the space of modular forms of level~1 to show that  of James and Ono~\cite{james-and-ono} to show that the other~$T^{1,1}_{l,k}$ are irreducible.

The situation in characteristic~$p$ is rather different. We see from Table~\ref{table-1} that the Hecke polynomials~$T^{1,1}_{l,k}$ factor into linear factors for all the primes less than~19. In this paper we will investigate such examples where a ``mod~$p$ version of Maeda's conjecture'' fail so completely.

We quote some results from~\cite{conrey-farmer-wallace}, which will show that the factorization of the Hecke polynomials~$T^{N,\chi}_{l,k}$ modulo~$p$ is periodic. Firstly, we will show that certain mod~$p$ reductions of Hecke polynomials divide one another. Let~$B_k(N,\chi)$ be a basis of~$S_k(\Gamma_0(N),\chi)$.
\begin{lemma}[Lemma 1,~\cite{conrey-farmer-wallace}]
\label{lemma-cfw}
Let~$[T_n]_{k,N,\chi}$ be the matrix of the Hecke operator~$T_n$ with respect to the basis~$B_k(N,\chi)$, and let~$q=p-1$ if~$p \ne 2$ and~$2$ if~$p = 2$. Then~$[T_n]_{k,N,\chi}$ is block upper-triangular and
\[
[T_n]_{k,N,\chi} \subset [T_n]_{k+q,N,\chi},
\]
where the smaller matrix appears as a block in the upper left corner of the larger matrix. In particular, we have
\[
T^{N,\chi}_{n+q} \equiv g(x) \cdot T^{N,\chi}_{n,k}\mod p
\]
for some polynomial~$g(x)$.
\end{lemma}
\begin{proof}
If~$p \ge 5$, we can follow the argument of~\cite{conrey-farmer-wallace}; recall that the classical Eisenstein series~$E_{p-1}$ has Fourier expansion at~$\infty$ which is congruent to~$1$ modulo~$p$, so there is a an inclusion
\[
S_k(\Gamma_0(N),\chi) \subset S_{k+p-1}(\Gamma_0(N),\chi)
\]
given by multiplication by~$E_{p-1}$. We can choose~$B_k(N,\chi)$ to respect this inclusion, and we also note that the Hecke operator~$T_n \mod p$ respects this inclusion map because using the definition of the Hecke operators we have
\[
(T_n f)(z) = \frac{1}{n}\sum_{ad=n} \chi(a) a^k \sum_{0 \le b < d} f \left(\frac{az+b}{d}\right)
\]
and~$a^k \equiv a^{k+p-1} \mod p$.

We note that if~$p=2$ or~$p=3$, then the argument involving~$E_{p-1}$, which is the argument at the bottom of page~124 of~\cite{conrey-farmer-wallace} is not applicable, as there is no level~1 modular form of weight~$1$ or~$2$. The solution here is to consider the spaces of modular forms as being mod~$p$ modular forms in the sense of Katz~\cite{katz}; this then means that the Hasse invariant~$A$ is a nonzero weight~$p-1$ mod~$p$ modular form, and then there is a map from weight~$k$ to weight~$k+q$ given by multiplication by~$A$ (if~$p =3$) or~$A^2$ (if~$p=2$). The rest of the proof does not depend on~$p \ge 5$.
\end{proof}
Using Lemma~\ref{lemma-cfw}, we see that there is a sequence of polynomials~$f_j \in (\Z/p\Z)[x]$ depending only on~$N$, $\chi$, $p$, $l$ and~$(k \mod p-1)$ such that
\begin{equation}
\label{the-fi}
T_{l,k}^{N,\chi}(x) \equiv \prod_{j=1}^J f_j(x) \mod p.
\end{equation}

%Explain the~$f_i$ here, say why they are periodic.
%\begin{equation}
%\label{the-fi-old}
%\prod_{i=1}^n f_i
%\end{equation}
It is in fact the case that a finite amount of computation suffices to compute any of the~$T_{l,k}^{N,\chi}(x)$, by the following proposition:
\begin{proposition}[Proposition~1, \cite{conrey-farmer-wallace}]
The sequence of polynomials~$f_i$ in~\eqref{the-fi} is periodic.
\end{proposition}
\begin{proof}
See~\cite{conrey-farmer-wallace}.
\end{proof}
One elementary observation which follows from~\eqref{the-fi} is that the degrees of the~$f_j$ are bounded above by
\begin{equation}
\label{bound-eqn}
M:=\max_k \left(\dim S_{k+p-1}(\Gamma_0(N),\chi)-\dim S_{k}(\Gamma_0(N),\chi)\right);
\end{equation}
if we have that~$p = 2$, then this formula is still true, but we see that this formula is not optimal, because one of these two spaces is zero. We can obtain a better bound by considering
\begin{equation}
M:=\max_k \left(\dim S_{k+2}(\Gamma_0(N),\chi)-\dim S_{k}(\Gamma_0(N),\chi)\right)
\end{equation}
or (if one is willing to consider mod~$2$ modular forms)
\begin{equation}
M:=\max_k \left(\dim S_{k+1}(\Gamma_0(N),\chi;\F_2)-\dim S_{k}(\Gamma_0(N),\chi;\F_2)\right)
\end{equation}
instead.

In Section~3 of~\cite{conrey-farmer-wallace}, it is shown that the~$T^{1,1}_{l,k}$ factor completely for~$l \le 7$ and~$l = 13$; this follows from~\eqref{bound-eqn}, because~$M=1$ in these particular cases. In Section~\ref{table-section}, we will extend their results to higher levels.

A natural question that arises is:
\begin{question}
Are there finitely many prime levels~$N$ and primes~$p$ such that the polynomials~$T_{l,k}^{N,1} \mod p$ split completely? If there are only finitely many, is the table in Section~\ref{table-section} complete?
\end{question}
In Section~\ref{heuristics-section}, we will give a heuristic argument suggesting that there should indeed be only finitely many levels~$N$ and primes~$p$ such that~$T_{l,k}^{N,1} \mod p$ splits completely.

A natural extension to this question is to consider the situation for other finite fields~$\F_{p^r}$; it is clear that if we enlarge the field of definition, then we will find more Hecke polynomials which split completely, so we have to specify a fixed~$\F_q$ for the question to be meaningful.
\section{Tables of results}
\label{table-section}
In Table~\ref{table-1}, we see a list of all prime levels~$N$ (with~$p\nmid N$), together with~$N=1$, less than~1000 where all of the~$T^{N,1}_{l,k}$ factor completely modulo~$p$, for~$p$ less than~41.
%\vskip 0.2in
\begin{table}[ht!]
\begin{tabular}{|c|c|}
\hline
$N$ & primes $p$ where the~$T^{N,1}_{l,k}$ factor completely\\
\hline
1 & 2,\;3,\;5,\;7,\;11,\;13,\;17,\;19 \\
\hline
2 & 2,\;3,\;5,\;7\\
\hline
3 & 2,\;5,\;7,\;11\\
\hline
5 & 2,\;3 \\
\hline
7 & 2,\;3 \\
\hline
11 & 2,\;3 \\
\hline
13 & 2 \\
\hline
17 & 2 \\
\hline
19 & 2,\;3 \\
\hline
29 & 2 \\
\hline
41 & 2 \\
\hline
53 & 2 \\
\hline
89 & 2 \\
\hline
\end{tabular}
\caption{Levels~$N$ where all of the~$T^{N,1}_{l,k}$ factor completely}
\label{table-1}
\end{table}
%\vskip0.2in
\noindent This table was computed using {\sc Magma}~\cite{magma}.

In~\cite{conrey-farmer-wallace}, it is shown that if the level is~1 and~$p \in \{2,3,5,7,11 \}$ then the~$T^{1,1}_{l,k}$ will factor completely modulo~$p$, because each of the~$f_i$ will have degree~1. We notice that for our largest example, level~89 and~$p=2$, the dimension of~$S_2(\Gamma_0(89))$ is~7, so therefore it is possible that quite large degree Hecke polynomials can still split modulo~$p$; this, again, differs sharply from the characteristic~0 case, where it is believed that the only eigenforms for~$\sltz$ which have coefficients in~$\Q$ are those which are forced to be defined over~$\Q$ by dimension considerations, such as the~$\Delta$-function.

Table~\ref{table-2} is the corresponding table showing those prime levels~$N$ less than~1000, together with~$N=4$, where all of the~$T^{N,\chi}_{l,k}$ split completely modulo~$p$ for~$\chi=\left(\frac{\cdot}{N}\right)$.
\begin{table}[ht!]
\begin{tabular}{|c|c|}
\hline
$N$ & primes~$p$ where all the~$T^{N,\chi}_{l,k}$ factor completely\\
\hline
3 & 2,5,7 \\
\hline
4 & 3,5 \\
\hline
5 & 2,3 \\
\hline
7 & 2,3,5 \\
\hline
11 & 2,3,5 \\
\hline
13 & 2,3 \\
\hline
17 & 2\\
\hline
19 & 2 \\
\hline
23 & 3 \\
\hline
29 & 2 \\
\hline
31 & 3 \\
\hline
41 & 2 \\
\hline
53 & 2 \\
\hline
\end{tabular}
\label{table-2}
\caption{Levels~$N$ where all of the~$T^{N,\chi}_{l,k}$ factor completely}
\end{table}
We note that for the largest of these levels, the dimension of the space of modular forms~$S_2(\Gamma_0(53),\left(\frac{\cdot}{N}\right))$ is~8, which is very close to the largest dimension, 7, that we found above. We also note that there are two levels, 23 and~31, where the Hecke polynomials factor completely modulo~3 but not modulo~2.

\section{Heuristics}
\label{heuristics-section}
We will briefly outline elementary arguments for why there should be only finitely many spaces of modular forms~$S_k(\Gamma_0(N))$ where all the Hecke polynomials factor completely over~$\F_p$, for~$N$ prime. Let us choose such an~$N$.

It is well-known (see~\cite{cohen-oesterle}, for instance) that the dimension of these spaces of modular forms tends to~$\infty$ as~$N$ tends to~$\infty$ (in a rather irregular way; see Remark~6.3 of~\cite{stein-computational-book}), and that all but finitely many spaces~$S_k(\Gamma_0(N))$ are greater than any given integer. This means that the quantity~$M$ in~\eqref{bound-eqn} is increasing, so the spaces which have Hecke polynomials which factor automatically by dimension considerations are finite in number.

However, we noted in the previous section that there are examples where the Hecke polynomials factor completely when~$M \ne 1$; we note that in our tables these become rarer and rarer as the dimension becomes higher.

Let us assume that the Hecke polynomials~$T_{2,k}^{N,1} \mod p$ are randomly distributed amongst polynomials of their degree~$d$. There are~$p^d$ monic polynomials of degree~$d$, and the number of polynomials of degree~$d$ which split completely is polynomial in~$d$; for instance, if~$p=2$, it is~$d+1$, and if~$p=3$, it is~$(d+1)(d+2)/2$. Therefore, as~$N$ increases, the degree of the first polynomial~$f_1$ in the expansion~\eqref{the-fi} increases, and this makes it more likely that the polynomial~$f_1$ does not factor completely.

%If~$f_1$ does not factor completely, then this means that none of the Hecke polynomials will factor completely. Even if it does factor completely, then it is still possible that one of the later~$f_i$ will fail to factor completely.

\section{Acknowledgements}
The computations contained in this paper were performed using {\sc Magma}~\cite{magma} on the machine Meccah owned by William Stein. The author is grateful for being given access to this machine.

The author would like to acknowledge helpful conversations with Edray Goins and David Farmer.
%\bibliographystyle{plain}
%\bibliography{mynewbib}

\end{document}